\documentclass[a4paper]{amsart}

\usepackage{amsmath,amsthm,amssymb,enumerate,graphicx,soul, overpic,enumitem}
\usepackage[usenames,dvipsnames,svgnames]{xcolor}
\usepackage[margin=1.3in]{geometry}
\usepackage{hyperref} 
\hypersetup{colorlinks=true,
allcolors=purple} 

\title{Orthosystoles and orthokissing numbers}
\author{Ara Basmajian} 
\address{(Ara Basmajian) The Graduate Center, CUNY, 365 Fifth  Ave., N.Y.,N.Y., 10016,
USA and Hunter College, CUNY, 695 Park Ave., N.Y.,N.Y., 10065, USA}
\thanks{}
\email{abasmajian@gc.cuny.edu}
\author{Federica Fanoni}
\address{(Federica Fanoni) CNRS, Univ Paris Est Creteil, Univ Gustave Eiffel, LAMA,
UMR8050, F-94010 Creteil, France}
\thanks{}
\email{federica.fanoni@u-pec.fr}
\date{\today}

\DeclareMathOperator{\sys}{sys}
\DeclareMathOperator{\osys}{osys}

\DeclareMathOperator{\teich}{Teich}
\DeclareMathOperator{\kiss}{kiss}
\DeclareMathOperator{\okiss}{okiss}
\DeclareMathOperator{\dteich}{d_T}

\DeclareMathOperator{\ms}{\mathcal{M}}
\DeclareMathOperator{\inj}{inj}

\newcommand{\R}{\mathbb{R}}

\newcommand{\A}{\mathcal{A}}
\renewcommand{\st}{\;|\;}
\newcommand{\ssm}{\smallsetminus}

\newtheorem{thm}{Theorem}[section]
\newtheorem{lemma}[thm]{Lemma}
\newtheorem{prop}[thm]{Proposition}

\newtheorem{cor}[thm]{Corollary}

\newtheorem{thmintro}{Theorem}

\theoremstyle{definition}
\newtheorem{rmk}[thm]{Remark}


\begin{document}
\begin{abstract}
For hyperbolic surfaces with geodesic boundary, we study the \emph{orthosystole}, ie the length of a shortest essential arc from the boundary to the boundary. We completely characterize local and global maxima for the orthosystole function among surfaces of given total boundary length, recovering  and extending  work by Bavard. For surfaces with fixed individual boundary lengths, we construct surfaces with large orthosystole and show that their orthosystole grows, as the genus goes to infinity, at the same rate as Bavard's upper bound.
\end{abstract}

\maketitle

\section{Introduction}

Systoles -- shortest closed geodesics -- of closed hyperbolic surfaces have been widely studied and have been fundamental in understanding hyperbolic surfaces and their moduli spaces. A simple area argument provides an upper bound on the systole length which is asymptotic to \(2\log g\), as the genus \(g\) goes to infinity. Surprisingly, this naive bound provides the correct order of growth, as various authors have constructed sequences of closed hyperbolic surfaces with systole length growing logarithmically in the genus (\cite{brooks_injectivity}, \cite{bs_period}, \cite{ksv_logarithmic}, \cite{petri_hyperbolic}, \cite{pw_graphs}, \cite{lp_random}). The best constructions (\cite{bs_period}, \cite{ksv_logarithmic}) have systole length growing at least as \(\frac{4}{3}\log g\). It is a well known open problem to understand the gap between the best construction and the best known upper bound (see also \cite{bavard_disques} and \cite{fbp_linear} for improvements on the naive bound).

In analogy with sphere packing problems in Euclidean space, Schmutz Schaller introduced and studied kissing numbers -- numbers of systoles -- of hyperbolic surfaces. Also for kissing numbers, upper bounds (\cite{parlier_kissing}, \cite{fp_systoles}, \cite{fbp_linear}) and constructions of surfaces with large kissing numbers (\cite{schmutz_arithmetic}, \cite{schmutz_compact}, \cite{schmutz_extremal}, \cite{brooks_platonic}, \cite{bmp_systole}) are known. Furthermore, it is known (\cite{schmutz_congruence}) that surfaces which are local maxima for the systole  must have large kissing number as well.

If we consider compact hyperbolic surfaces with boundary, it is natural to consider another collection of geodesics, instead of closed ones: those starting and ending at the boundary components. Given a hyperbolic surface \(X\) with geodesic boundary, we define the \emph{orthosystole} to be the length of a shortest geodesic from boundary to boundary, and we denote it by \(\osys(X)\). Moreover, we call \emph{orthokissing number} (and we denote it by \(\okiss(X)\)) the number of geodesics from boundary to boundary having minimal length.

The first observation is that there are no interesting upper and lower bounds for \(\osys\) depending only on the topology of the surface: if \(S\) has negative Euler characteristic, the orthosystole can be arbitrarily small and arbitrarily large (see Lemma \ref{lem:no-bounds}). As a consequence, we will discuss bounds on the orthosystole for surfaces with restrictions on the boundary lengths. Two natural possibilities are to fix each boundary length or to fix their sum, and it turns out that the two situations are quite different. 

If we fix the sum of the boundary lengths, a sharp upper bound has been proven by Bavard in \cite{bavard_anneaux}:

\begin{thm}[Bavard]\label{thm:Bavard}
Let S  be a  surface of signature \((g,n)\)  with 
$\chi (S) <0$.  For every $X \in \ms(S)$  having  total boundary length \(L\)  
\[
\osys(X)
\leq 2\sinh^{-1}\left(\frac{1}{2\sinh\left(\frac{L}{24g-24+12n}\right)}\right)
\]
with equality if and only if  the orthosystoles decompose $X$ into hyperbolic right-angled hexagons. Moreover, for every 
such \((g,n)\) and  every \(L>0\) equality is attained.
\end{thm}

Note that Bavard's work is more general and is phrased quite differently. The theorem just stated corresponds to a special case of parts (1), (3) and (4) of \cite[Th\'eor\`eme 1]  {bavard_anneaux}. We refer to the Appendix for an explanation of the correspondence between Bavard's formulation and the theorem above.

Our first result extends Bavard's work  by showing that global and local maxima for the orthosystole function coincide.

We denote by \(\ms(S;L)\) the space of hyperbolic structures on \(S\) with fixed total boundary length \(L\).

\begin{thmintro}\label{thm:n=1}
Let \(S\) be a surface of signature \((g,n)\) with 
$\chi(X) <0$ and let \(L>0\). For \(X\in\ms(S;L)\), the following are equivalent:
\begin{enumerate}
\item \(X\) is a global maximum for the orthosystole function on \(\ms(S;L)\);
\item \(X\) is a local maximum for the orthosystole function  on \(\ms(S;L)\);
\item \(\okiss(X)=6g-6+3n\);
\item \(\osys(X) = 2\sinh^{-1}\left(\frac{1}{2\sinh\left(\frac{L}{24g-24+12n}\right)}\right)\).
\end{enumerate}
Moreover,  for $n=1$ and  any  $g \geq 1$ the number of local (or global) maxima in each moduli space is   exactly
\[\frac{2 (6g-5)!}{12^g g!(3g-3)!}.\]
\end{thmintro}

We remark that our proof is very different from the one of Bavard ---  in particular, we do not rely on his work to prove Theorem \ref{thm:n=1}. Note  also that condition (3) in Theorem
\ref{thm:n=1} is equivalent to the orthosystoles decomposing  the surface into hexagons (see Section  \ref{sec:hexdec}).

 For the analogous problem for the systole function, global maxima are not known (except in a few low complexity cases) and it is a hard problem even to construct local maxima (see \cite{fbr_local}).
 
 If, instead of fixing the sum,  we fix each boundary length,  we can not always find  a hyperbolic surface with 
 \(\okiss(X)=6g-6+3n\) (as discussed at the beginning of Section 
 \ref{sec:multiplecomponents}). So Bavard's result and Theorem 
 \ref{thm:n=1} are not always sharp.
 Still, while we cannot completely describe global maxima, we are able to show that the orthosystole admits a global maximum,   construct examples of surfaces with large orthosystole, and show that local maxima of the orthosystole function  have relatively large orthokissing number.  

We denote by \(\ms(S;\ell_1,\dots,\ell_n)\) the space of  hyperbolic structures on $S$ with boundary lengths 
$\ell_1,\dots,\ell_n$.

\begin{thmintro}\label{thm:natleast2}
Let \(S\) be a surface of signature \((g,n)\) and negative Euler characteristic. Fix \(0<\ell_1\leq \dots\leq\ell_n\).
\begin{enumerate}
\item The function \(\osys:\ms(S;\ell_1,\dots,\ell_n)\to \R\) admits a maximum.
\item Suppose \(g\geq n\geq 2\). Then there is a surface \(X\in \ms(S;\ell_1,\dots,\ell_n)\) with
\[\osys(X)\geq \cosh^{-1}\left(\cfrac{\cosh\left(\cfrac{\ell_n}{6\left\lfloor \frac{g}{n}\right\rfloor}\right)}{\cosh\left(\cfrac{\ell_n}{6\left\lfloor \frac{g}{n}\right\rfloor}\right)-1}\right).\]
\item If \(X\in \ms(S;\ell_1,\dots,\ell_n)\) is a local maximum for 
\(\osys\), then \(\okiss(S)\geq 2g-2+n\).
\end{enumerate} 
\end{thmintro}

While the fact  that \(\osys\) admits a maximum might seem a triviality, it is more subtle than in the case of the systole. Indeed, continuity of the systole function and Mumford's compactness theorem readily imply the existence of global maxima for the systole function in each moduli space. On the other hand, we cannot replace the systole by the orthosystole in Mumford's criterion (see Lemma \ref{lem:noMumford}). So to prove (1) we will need to apply a length-expansion result (\cite{thurston_spine},\cite{parlier_lengths},\cite{pt_shortening}).

Note furthermore that if \(n\) and the \(\ell_i\) are fixed and we let the genus go to infinity, the surfaces we construct in part (2) have orthosystole  growing as \(2\log g\), which matches the asymptotics of the upper bound given in Theorem \ref{thm:Bavard}.
On the other hand (see Section \ref{sec:multiplecomponents}), these  surfaces  aren't even \emph{local} maximizers for the orthosystole. 

This article is structured as follows: after defining the objects we are interested in and recalling or proving some facts we will need (Section \ref{sec:prereq}), in Section \ref{sec:hexdec}  we  discuss maximal collections of pairwise disjoint and pairwise non-homotopic arcs, which will play a fundamental role in our work. Section \ref{sec:bounds} is dedicated to proving some basic facts about orthosystoles and orthokissing numbers and part (1) of Theorem \ref{thm:natleast2} (see Corollary \ref{cor:maxisrealized}). Section \ref{sec:oneboundary} is concerned with the case of fixed total boundary length and the proof of Theorem \ref{thm:n=1}. The case of fixed  boundary lengths is studied in  Section \ref{sec:multiplecomponents}, where part (2) of Theorem \ref{thm:natleast2} is proven. Finally, in Section \ref{sec:osys&okiss} we show part (3) of Theorem \ref{thm:natleast2}.

\section*{Acknowledgements}
The  authors  would like to thank Hugo Parlier and Bram Petri for useful discussions.  We are grateful to the referee for the insightful comments.

A.B.\ was partially supported by PSC CUNY Award 65245-00 53 and the Simons Foundation (359956, A.B.).

\section{Prerequisites}\label{sec:prereq}
Throughout the rest of this manuscript we assume 
our surfaces to be compact with non-empty boundary and negative Euler characteristic. Hyperbolic surfaces are assumed metrically complete with geodesic boundary.

Given a compact surface \(S\) of signature \((g,n)\), we denote by \(\teich(S)\) the Teichm\"uller space of hyperbolic surfaces with geodesic boundary homeomorphic to $S$, where the boundary length is allowed to vary. Given \(L>0\), \(\teich(S; L)\) is the subspace of surfaces whose total boundary length is \(L\). For \(\ell_1,\dots,\ell_n>0\), \(\teich(S;\ell_1,\dots,\ell_n)\) is the subspace of surfaces whose \(i\)-th boundary component has length \(\ell_i\), for \(i=1,\dots,n\). We denote similarly the moduli space  \(\ms(S)\) and its subspaces  \(\ms(S; L)\) and \(\ms(S;\ell_1,\dots,\ell_n)\). Moreover, \(\ms_\varepsilon(S;\ell_1,\dots,\ell_n)\) denotes the subset of \(\ms(S;\ell_1,\dots,\ell_n)\) given by surfaces all of whose non-boundary parallel closed geodesics have length at least \(\varepsilon\) -- the so-called \emph{(\(\varepsilon\)-)thick part} of moduli space. Using a doubling argument and  the Collar Lemma (see Lemma 
\ref{lem:collarlemma}), a straightforward application of  Mumford's compactness theorem \cite{mumford_remark} 
 shows  that \(\ms_\varepsilon(S;\ell_1,\dots,\ell_n)\)   is compact for every $\varepsilon >0$.

Let \(X\in\teich(S)\). Its \emph{systole} \(\sys(X)\) is the length of a shortest closed geodesic which is not a boundary component. Its \emph{orthosystole} \(\osys(X)\) is the length of a shortest orthogeodesic (a geodesic segment between boundary components which intersects the  boundary at right-angles).  We will regularly abuse notation and use the terms systole and orthosystole also for a simple closed geodesic (not boundary parallel) or an orthogeodesic of minimal length. The  \emph{kissing number} of $X$,  \(\kiss(X)\),  is the number of systoles and the  \emph{orthokissing number} of $X$, \(\okiss(X)\),  is the number of  orthosystoles. 

Given a homotopically nontrivial simple closed curve \(\gamma\), we denote by \(\ell_{\gamma}(X)\) the length of the unique shortest closed curve in the free homotopy class of \(\gamma\), ie of the unique simple closed geodesic in the class.

Denote by \(\A=\A(S)\) the collection of essential arcs in \(S\) from boundary to boundary, up to homotopy relative to the boundary. For any \(\alpha\in\A\), we denote by \(\ell_{\alpha}(X)\) the length of the unique shortest arc  in the class, ie of the unique orthogeodesic in the class.

\begin{lemma}\label{lem:hextrig}
Let \(a, \gamma, b, \alpha, c, \beta\) be the lengths of the sides of a right-angled hyperbolic hexagon listed  in counter clockwise order, and let \(h_a\) be the length of the orthogonal between
the sides of length $a$ and $\alpha$.  Let \(A=\cosh(a),B=\cosh(b)\) and \(C=\cosh(c)\) and set 
\begin{align*}s(A,B,C):=&\cosh^{-1}\left(\frac{A+BC}{\sqrt{(B^2-1)(C^2-1)}}\right)+\\
&+\cosh^{-1}\left(\frac{B+AC}{\sqrt{(A^2-1)(C^2-1)}}\right)+\cosh^{-1}\left(\frac{C+AB}{\sqrt{(A^2-1)(B^2-1)}}\right).
\end{align*}
Then:
\begin{enumerate}[itemsep=5pt]
\item \(\displaystyle\alpha=\cosh^{-1}\left(\frac{A+BC}{\sqrt{(B^2-1)(C^2-1)}}\right)\);
\item \(\alpha + \beta + \gamma=s(A,B,C)\);
\item \(\displaystyle \frac{\partial s}{\partial A}(A,B,C)=\frac{A-1-B-C}{(A-1)\sqrt{A^2+B^2+C^2+2ABC-1}}\);
\item For every $B,C>1$, \(\displaystyle \lim_{A\to 1}s(A,B,C)=+\infty=\lim_{A\to +\infty}s(A,B,C)\);
\item \(\displaystyle \cosh(h_a)=\frac{\sqrt{A^2+B^2+C^2+2ABC-1}}{\sqrt{A^2-1}}\); if \(b=c\),
\(\displaystyle \sinh(h_a)=\frac{\cosh(b)}{\sinh(a/2)}\).
\end{enumerate}
\end{lemma}

\begin{proof}
Standard  formulas   (see for instance \cite[Chapter 2]{buser_geometry}) give us the lengths of  \(\alpha, \beta,\) and  \(\gamma\) and hence their sum  \(s(A,B,C)\). The formula for the derivative of \(s\) follows by explicit computation.

If \(A\to 1\), 
\[s(A,B,C)\geq \cosh^{-1}\left(\frac{C+AB}{\sqrt{(B^2-1)(A^2-1)}}\right)\to \infty.\]  If \(A\to\infty\)
\[s(A,B,C)\geq \cosh^{-1}\left(\frac{A+BC}{\sqrt{(B^2-1)(C^2-1)}}\right)\to \infty.\]

The expressions for \(h_a\) follow again from standard hyperbolic trigonometry (\cite[Chapter 2]{buser_geometry}), by looking at the two right-angled pentagons obtained by cutting the hexagon along the orthogonal from the side of length \(a\) to the opposite side.
\end{proof}

\begin{figure}[h]
\begin{center}
\begin{overpic}[scale=1.5]{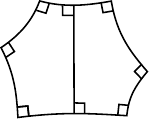}
\put(47,-5){\(a\)}
\put(86,55){\(b\)}
\put(8,60){\(c\)}
\put(47,80){\(\alpha\)}
\put(0,20){\(\beta\)}
\put(92,12){\(\gamma\)}
\put(53,40){\(h_a\)}
\end{overpic}
\caption{A right-angled hexagon}
\end{center}
\end{figure}
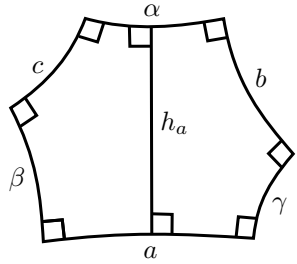

We will also use the Collar Lemma (\cite{keen_collars}; see also \cite[Chapter 4]{buser_geometry}):
\begin{lemma}[Collar Lemma]\label{lem:collarlemma}
A simple closed geodesic \(\alpha\) in a hyperbolic surface \(X\) has an embedded collar of width
\[w(\alpha)=\sinh^{-1}\left(\cfrac{1}{\sinh\left(\frac{\ell_\alpha(X)}{2}\right)}\right).\]
In particular any simple closed geodesic or orthogeodesic transversely intersecting \(\alpha\) has length at least \(w(\alpha)\) (\( 2w(\alpha)\) if \(\alpha\) is not a boundary component).
\end{lemma}

The following is due to Wolpert  in the closed surface case (\cite{wolpert_length}) and deduced in   \cite[Theorem 2.1]{lpst_length}) for the case with boundary:
\begin{thm}\label{thm:Kqc}
If \(f:X\to Y\) is a \(K\)-quasiconformal homeomorphism between compact hyperbolic surfaces with boundary and \(\gamma\subset X\) is either an essential simple closed curve or an essential arc from boundary to boundary, then
\[\frac{1}{K} \leq \frac{\ell_{f(\gamma)}(Y)}{\ell_{\gamma}(X)}\leq K.\]
\end{thm}

We will also use a length-expansion result (see \cite{thurston_spine}, \cite{parlier_lengths} and \cite{pt_shortening}):
\begin{thm}\label{thm:lengthexpansion}
Let \(X\) be a hyperbolic surface of signature \((g,n)\), and let \(\gamma_1,\dots,\gamma_n\) be the boundary geodesics of \(X\). For every \((\varepsilon_1,\dots,\varepsilon_n)\in \R^n \) where 
\(\varepsilon_i \geq 0 \) for all \(i\), and at least one \(\varepsilon_i \)
is not zero, 
there exists a hyperbolic surface \(Y\) homeomorphic to $X$ with boundary geodesics of lengths \(\ell_{\gamma_1}(X)+\varepsilon_1,\dots, \ell_{\gamma_n}(X)+\varepsilon_n\) such that each  closed geodesic  in \(Y\) is  strictly longer  than the corresponding geodesic  in  \(X\).
\end{thm}

\begin{rmk}\label{rmk:lengthexpansion}
For surfaces \(X,Y\) as in Theorem \ref{thm:lengthexpansion}, if \(\varepsilon_i=\varepsilon_j=0\) and \(\alpha\) is an orthogeodesic from \(\gamma_i\) to \(\gamma_j\), \(\ell_{\alpha}(Y)>\ell_{\alpha}(X)\). This follows by the proof of the theorem given in \cite{parlier_lengths}.
\end{rmk}

\section{Hexagon decompositions}\label{sec:hexdec}

A \emph{hexagon decomposition} \(H\) of a  surface \(S\)  is a maximal collection of pairwise non-homotopic and disjoint essential arcs. If \(X\in\teich(S)\), the orthogeodesic representatives of the arcs in \(H\) cut \(X\) into a union of right-angled hexagons. By Euler characteristic considerations, a hexagon decomposition of a surface of signature \((g,n)\) contains \(3|\chi(S)|=6g-6+3n\) arcs and cuts the surface into \(2|\chi(S)|=4g-4+2n\) (topological) disks.

Ushijima \cite{ushijima_canonical} showed:
\begin{thm}\label{thm:parametrization}
Let \(S\) be a surface of signature \((g,n)\) and \(H=\{\alpha_1,\dots\alpha_{6g-6+3n}\}\) a hexagon decomposition. Then
\begin{align*}
\varphi_H:\teich(S)&\to\R_{>0}^{6g-6+3n}\\
X&\mapsto (\ell_{\alpha_i}(X))_{i=1}^{6g-6+3n}
\end{align*}
is a homeomorphism. 
\end{thm}

\begin{rmk} Note that hexagon decompositions in \cite{ushijima_canonical}
are called truncated triangles. 
\end{rmk}

For the sake of completeness, we provide here an alternative proof:
\begin{proof}
Since right-angled hexagons are determined by the lengths of three alternating sides, and any three positive lengths determine a right-angled hexagon, the map is a bijection. Theorem \ref{thm:Kqc} implies that \(\varphi_H\) is continuous. Conversely, Bishop \cite{bishop_quasiconformal} shows that given two right-angled hexagons of alternating side lengths \(a,b,c\) and \(a',b',c'\), there is a quasiconformal homeomorphism from one to the other, sending the side of length \(a\) (respectively \(b,c\)) to the side of length \(a'\) (respectively, \(b',c'\)), linear on the sides, whose quasiconformal constant depends on \(a,b,c,\max\{|a-a'|,| b-b'|,|c-c'|\}\) and goes to one as \(\max\{|a-a'|,| b-b'|,|c-c'|\}\) goes to zero. By using these maps on each hexagon determined by \(H\), we get a quasiconformal map \(\varphi_H^{-1}(x_1,\dots,x_{6g-6+3n})\to \varphi_H^{-1}(x'_1,\dots,x'_{6g-6+3n})\), whose quasiconformal constant depends on \(x_1,\dots,x_{6g-6+3n},\max\{|x_i-x'_i|\st i=1,\dots, 6g-6+3n\}\) and goes to one as \(\max\{|x_i-x'_i|\st i=1,\dots, 6g-6+3n\}\) goes to zero. This proves the continuity of \(\varphi_H^{-1}\).
\end{proof}

Note that this means in particular that for every arc \(\alpha\in\A\), the function
\begin{align*}
\teich(S)&\to\R\\
X&\mapsto \ell_{\alpha}(X)
\end{align*}
is continuous.

Given a surface \(S\) of signature \((g,n)\)  and a hexagon decomposition \(H=\{\alpha_1,\dots,\alpha_{6g-g+3n}\}\), denote by \(\partial_H\) the collection of triples \(\{i,j,k\}\) (where two indices might coincide) of indices such that \(\alpha_i,\alpha_j,\alpha_k\in H\) are three sides of some hexagon defined by \(H\).  Define the function\footnote{Here we are slighlty abusing notation, since \(\{i,j,k\}\) is an unordered triple, so \(s(\cosh(x_i),\cosh(x_j),\cosh(x_k))\) is technically not defined, but since \(s\) is symmetric, we can choose any ordering for \(i,j,k\) to compute \(s\).}
\begin{align*}
F_H:\R^{6g-6+3n}&\to \R\\
x&\mapsto \sum_{\{i,j,k\}\in\partial_H} s(\cosh(x_i),\cosh(x_j),\cosh(x_k)).
\end{align*}

We have

\begin{lemma}\label{lem:image}
Let \(S\) be a surface of signature \((g,n)\) and \(H=\{\alpha_1,\dots\alpha_{6g-6+3n}\}\) a hexagon decomposition. For every \(L>0\)
\[\varphi_H(\teich(S; L))=\left\{x\in\R_{>0}^{6g-6+3n}\;\middle|\;~ F_H(x)=L\right\}.\]
\end{lemma}

\begin{proof}
This is a consequence of Lemma \ref{lem:hextrig} and the fact that the boundary of \(S\) is the union of all the sides not belonging to \(H\) in the hexagons in the decomposition given by \(H\).
\end{proof}

\section{Orthosystoles and orthokissing number: general bounds}\label{sec:bounds}

This section is dedicated to proving some basic properties of orthosystoles and part (1) of Theorem \ref{thm:natleast2}. We start by showing that there are no bounds on the orthosystole depending only on the geometry of the surface:

\begin{lemma}\label{lem:no-bounds}
Let \(S\) be a compact surface with $\chi(S)<0$ and nonempty boundary. Then
\[\inf_{X\in\ms(S)}\osys(X)=0 \;\;\;\text{and}\;\;\;\sup_{X\in\ms(S)}\osys(X)=\infty.\]
\end{lemma}
\begin{proof}
Given a pair of pants of boundary lengths \(a,b,b\), the length of the shortest orthogeodesic from and to the boundary component of length \(a\) goes to zero if \(a\to\infty\) and \(b\to 0\). By either gluing the two boundary components of length \(b\) to each other or attaching a subsurface to them, we get surfaces in each moduli space with arbitrarily short orthosystole. The fact that the supremum is infinite follows from the Collar Lemma, by letting all the boundary lengths go to zero.
\end{proof}

Systoles in a compact hyperbolic surface (which is not a pair of pants)  are simple.  Similarly,    orthosystoles are simple. On the other hand, while two systoles can intersect, we prove that orthosystoles are disjoint, which will imply a sharp  bound on the orthokissing number.

\begin{lemma}\label{lem:simple&disjoint}
Orthosystoles in a compact hyperbolic surface are simple and pairwise disjoint. In particular, a surface \(X\) of signature \((g,n)\) satisfies \(\okiss(X)\leq 6g-6+3n\).
\end{lemma}
\begin{proof}
Simplicity and disjointness follow from a standard surgery argument: if an orthosystole is not simple or two orthosystoles intersect, we can construct a shorter orthogeodesic using subarcs determined by some (self-)intersection point. This uses the fact that  a hyperbolic triangle can not have two right angles, and hence the geodesic subarcs  only intersect if they bound non-trivial topology. 

As there are at most \(6g-6+3n\) pairwise disjoint orthogeodesics on a surface of signature \((g,n)\), the bound on \(\okiss(X)\) follows immediately.
\end{proof}

The following criterion to detect orthosystoles will be crucial in the proofs of Theorems \ref{thm:n=1} and \ref{thm:natleast2}.

\begin{lemma}\label{lem:decosys}
Let \(X\) be a hyperbolic surface with a hexagon decomposition of orthogeodesics of the same length. Then the orthogeodesics in the hexagon decomposition are precisely the orthosystoles.
\end{lemma}
\begin{proof}
Denote by \(a\) the length of the arcs in the hexagon decomposition and let \(\alpha\) be the length of the other sides of a hexagon with alternating sides of length \(a\).

Let \(\gamma\) be an orthogeodesic of \(X\) not in the hexagon decomposition. It starts at some point \(p\in\partial X\), which needs to be contained in the interior of a side of length \(\alpha\) in some hexagon, and it is split into arcs \(a_1,\dots, a_k\) by the orthogeodesics in the hexagon decomposition, where \(k\geq 2\) (since \(\gamma\) is not in the hexagon decomposition). As there are no hyperbolic triangles with two right angles, \(a_1\) and \(a_k\) need to join opposite sides of a hexagon. In particular they both have length at least \(h\), where \(h\) is the orthogonal between opposite sides of a hexagon, so \(\ell(\gamma)\geq 2h\). Moreover by Lemma \ref{lem:hextrig}
\[\sinh(h)=\frac{\cosh(a)}{\sinh(a/2)}=\frac{2\sinh^2(a/2)+1}{\sinh(a/2)}>\sinh(a/2)\]
ie \(2h>a\), proving that the orthogeodesics in the hexagon decomposition are the unique ones of minimal length.
\end{proof}

Next we show that the orthosystole function is locally the minimum of finitely many continuous functions, and thus continuous.

\begin{lemma}\label{lem:finitelymanylengths}
For every \(X\in \teich(S)\), there is an open neighborhood \(U\) of \(X\) so that if \(Y\in U\)
\[\osys(Y)=\min\{\ell_{\alpha}(Y)\st \alpha\in A\},\]
where
\[A=\{\alpha\in\A\st\ell_\alpha(X)=\osys(X)\}.\]
In particular, the orthosystole function
\begin{align*}
\osys:\teich(S)&\to\R\\
X&\mapsto \osys(X)
\end{align*}
is locally the minimum of the lengths of finitely many orthogeodesics, and hence is continuous.
\end{lemma}

\begin{proof}
Let \(X\in\teich(S)\) and \(K>1\). Let \(\dteich\) be the Teichm\"uller distance on \(\teich(S)\). Then for every \(Y\in\teich(S)\) such that \(\dteich(X,Y)\leq \frac{\log K}{2}\), we have \(\osys(Y)\leq K\osys(X)\) (by Theorem \ref{thm:Kqc}), so
\[\osys(Y)=\min\{\ell_{\alpha}(Y)\st\alpha\in\A:\ell_{\alpha}(Y)\leq K\osys(X)\}.\]
Again by Theorem \ref{thm:Kqc}
\[\{\alpha\in\A\st\ell_{\alpha}(Y)\leq K\osys(X)\}\subset\{\alpha\in\A\st\ell_{\alpha}(X)\leq K^2\osys(X)\}.\]

Set \(B:=\{\alpha\in\A\st\ell_{\alpha}(X)\leq K^2\osys(X)\}\) and note that this is finite (by discreteness). Moreover, for every \(Y\in\teich(S)\) such that \(\dteich(X,Y)\leq 
\frac{\log K}{2}\), we have
\[\osys(Y)=\min\{\ell_{\alpha}(Y)\st\alpha\in B\}.\]

Let \(A:=\{\alpha\in\A\st\ell_{\alpha}(X)=\osys(X)\}\); note that \(A\subseteq B\). Moreover, there is \(C>1\) so that, for every \(\alpha\in A\) and \(\beta\in B\ssm A\) we have
\[C^2\ell_\alpha(X)\leq \ell_\beta(X).\]

Choose \(K'<\min\{K,C\}\). Set \(U:=\{Y\in\teich(S)\st\dteich(X,Y)\leq \frac{\log K'}{2}\}\). By definition of \(K'\) and Theorem \ref{thm:Kqc}

\[\ell_\alpha(Y)\leq K'\ell_\alpha(X)< \frac{K'}{C^2}\ell_\beta(X)<\frac{1}{K'}\ell_\beta(X)\leq \ell_\beta(Y).\]
So \[\osys(Y)=\min_{\gamma\in B}\ell_{\gamma}(Y)=\min_{\gamma\in A}\ell_{\gamma}(Y),\]
as required.

As lengths of orthogeodesics are continuous functions, \(\osys\) is continuous as well. 
\end{proof}

In the case of simple closed geodesics, it follows directly from Mumford's compactness criterion and continuity of the systole function that there are global systole maximizers in moduli or Teichm\"uller space. One might wonder if one could show an analogue of Mumford's criterion using orthogeodesic lengths, ie proving that for every \(\varepsilon>0\), the set
\[\{X\in\ms(S;\ell_1,\dots,\ell_n)\st \osys(X)\geq \varepsilon\}\]
is a compact subset of moduli space. Unfortunately, this is not the case, as we will show in Section \ref{sec:multiplecomponents}. This is the reason why we need a different approach to show that the orthosystole function has a global maximum. To this end, we prove:

\begin{prop}\label{prop:thinpart}
For every \(S\) of signature \((g,n)\) and for every \(\ell_1,\dots,\ell_n>0\) there is \(\varepsilon>0\) such that if \(X\in\ms(S;\ell_1,\dots,\ell_n)\) satisfies \(\sys(X)<\varepsilon\), then there is \(Y\in\ms(S;\ell_1,\dots,\ell_n)\) with \(\sys(Y)\geq \varepsilon\) and \(\osys(Y)>\osys(X)\).
\end{prop}

\begin{proof} Denoting the total boundary length of $X$ by 
\(\ell(\partial X)\),
choose 
$$\varepsilon<\min\left\{\sinh^{-1}(1),2\sinh^{-1}\left( 2\sinh\left(\frac{\ell(\partial X)}{24g-24+12n}\right)\right)\right\}. 
$$

Suppose \(\sys(X)<\varepsilon\). Let \[C_\varepsilon=\{\gamma\in\mathcal{C}(S)\st\ell_{\gamma}(X)<\varepsilon\}\] 
where \(\mathcal{C}(S)\) is the set of all essential closed curves in \(S\). 
Note that if a geodesic crosses a curve  
 of length at most \(\varepsilon\), by Lemma \ref{lem:collarlemma} it has length at least \[2\sinh^{-1}\left(\frac{1}{\sinh(\varepsilon/2)}\right)> 2\sinh^{-1}\left(\frac{1}{2\sinh\left(\frac{\ell(\partial X)}{24g-24+12n}\right)}\right).\]
 In particular, the geodesic representatives of curves in \(C_\varepsilon\)
 are disjoint and 
  by Theorem \ref{thm:Bavard} any orthosystole is disjoint from \(C_\varepsilon\).

Let \(X'\) be the surface obtained by cutting \(X\) along \(C_\varepsilon\). Apply Theorem \ref{thm:lengthexpansion} and Remark \ref{rmk:lengthexpansion} to \(X'\) to get a surface \(Y'\) where all curves corresponding to curves in \(C_\varepsilon\) have length \(\varepsilon\) and such that all curves in \(Y'\) are longer than those in \(X'\) and all orthogeodesics of \(Y'\) joining two curves corresponding to boundary components of \(X\) are longer or equal to the corresponding in \(X'\). Glue the boundary components of \(Y'\) which corresponded to the same curve in \(X\) to get a surface \(Y\). Let \(A\subset\A\) be a finite set such that
\[\osys(X)=\min\{\ell_\alpha(X)\st \alpha\in A\}\;\;\mbox{and}\;\;\osys(Y)=\min\{\ell_\alpha(Y)\st \alpha\in A\}.\]
Since orthosystoles are disjoint from \(C_\varepsilon\), we can assume that \(A\) contains only arcs disjoint from \(C_\varepsilon\). By construction \(\sys(Y)=\varepsilon\) and for every \(\alpha\in A\),
\[\ell_\alpha(X)<\ell_\alpha(Y),\]
so \(\osys(X)<\osys(Y)\).
\end{proof}

An immediate consequence is the following:
\begin{cor}\label{cor:maxisrealized}
The function \(\osys\) restricted to \(\ms(S;\ell_1,\dots,\ell_n)\) admits a maximum.
\end{cor}

\begin{proof}
By Proposition \ref{prop:thinpart}, there is \(\varepsilon>0\) such that
\[\sup\{\osys(X)\st X\in\ms(S;\ell_1,\dots,\ell_n)\}=\sup\{\osys(X)\st X\in\ms_\varepsilon(S;\ell_1,\dots,\ell_n)\}.\] As \(\ms_\varepsilon(S;\ell_1,\dots,\ell_n)\) is compact and \(\osys\) is continuous on moduli space (since it is a mapping class group invariant continuous function on Teichm\"uller space), \(\osys\) admits a maximum
\end{proof}

\section{Surfaces with fixed total boundary length}\label{sec:oneboundary}

This section is dedicated to the case of surfaces with fixed total boundary length and the proof of Theorem \ref{thm:n=1}. The main reason why fixing the sum of the boundary length is an advantage over fixing each boundary length is the fact that in this case we can always construct surfaces with a hexagon decomposition of orthosystoles and we can explicitly compute their orthosystole. The main difficulty in the proof of Theorem \ref{thm:n=1} is showing that local maxima for the orthosystole function need to have a hexagon decomposition of orthosystoles.

We start by showing the existence of surfaces with a hexagon decomposition of orthosystoles.

\begin{lemma}\label{lem:hexdecsamelength}
Let \(S\) be a surface of signature \((g,n)\) and \(H\) a hexagon decomposition of \(S\). Let \(L>0\). There is a unique \(X\in\teich(S;L)\) such that all arcs in \(H\) have the same length and
\[\osys(X)=f(g,L):=2\sinh^{-1}\left(\frac{1}{2\sinh\left(\frac{L}{24g-24+12n}\right)}\right).\]
\end{lemma}

\begin{proof}
Uniqueness of the hyperbolic structure follows from Theorem \ref{thm:parametrization}, so we just need to check existence. 

So set
\[a=2\sinh^{-1}\left(\frac{1}{2\sinh\left(\frac{L}{24g-24+12n}\right)}\right).\]
Glue \(4g-4g+2n\) copies of a hyperbolic hexagon of alternating sides of length \(a\) according to the combinatorics of \(H\), along the sides of length \(a\).

In each hexagon, the unglued sides have the same length \(\alpha\) given by
\[\cosh(\alpha)=\frac{\cosh(a)+\cosh^2(a)}{\sinh^2(a)}=\frac{\cosh(a)+\cosh^2(a)}{\cosh^2(a)-1}=\frac{\cosh(a)}{\cosh(a)-1}=\cosh\left(\frac{L}{12g-12+6n}\right).\]
The boundary of \(X\) is given by \(3(4g-4+2n)\) \(\alpha\)-sides (since there are \(4g-2\) hexagons), ie
\[\ell(\partial X)=(12g-12+6n)\alpha=L\]
and hence \(X\in\teich(S;L)\).
By Lemma \ref{lem:decosys}, the arcs in \(H\) are precisely the orthosystoles, so the statement about \(\osys(X)\) follows.
\end{proof}

\begin{rmk}\label{rmk:hexdec}
The function \(f(g, L)\) is independent of the hexagon decomposition. In particular non-homeomorphic hexagon decompositions yield non-isometric surfaces with the same orthosystole (indeed, an isometry between surfaces with hexagon decompositions of orthogeodesics of minimal length needs to send the hexagon decomposition of a surface to the hexagon decomposition of the other surface, since the arcs in the hexagon decompositions are the only ones of minimal length).
\end{rmk}

To show that local maxima for the orthosystole function have a hexagon decomposition of orthosystole, we will need the following technical lemma.

\begin{lemma}\label{lem:derivatives}
Let \(X\in\teich(S)\), where \(S\) has signature \((g,n)\) and \(H=\{\alpha_1,\dots,\alpha_{6g-6+3n}\}\) a hexagon decomposition of \(S\).
\begin{enumerate}
\item If \(\ell_{\alpha_i}(X)=\min_{j}\ell_{\alpha_j}(X)\), then \(\frac{\partial F_H}{\partial x_i}(\varphi_H(X))<0\).
\item If \(\frac{\partial F_H}{\partial x_j}(\varphi_H(X))=0\), then \(\ell_{\alpha_j}(X)\) is a global minimum or an inflection point of the function
\begin{align*}
F_{H,X,i}:\R_+&\to\R\\
y&\mapsto F_H(\ell_{\alpha_1}(X),\dots, \ell_{\alpha_{j-1}}(X),y,\ell_{\alpha_{j+1}}(X),\dots \ell_{\alpha_{6g-6+3n}}(X)).
\end{align*}
\end{enumerate}
\end{lemma}

\begin{proof}
We have
\begin{align*}
F_H:\R^{6g-6+3n}&\to \R\\
x&\mapsto \sum_{\{i,j,k\}\in\partial_H} s(\cosh(x_i),\cosh(x_j),\cosh(x_k)).
\end{align*}
For any \(i\), \(x_i\) appears in exactly two terms in the sum, those corresponding to the two hegaxons containing the orthogeodesic \(\alpha_i\). Suppose the two terms correspond to triples \(\{i,j,k\}\) and \(\{i,l,m\}\). Then
\[\frac{\partial F_H}{\partial x_i}=\frac{\partial s(\cosh(x_i),\cosh(x_j),\cosh(x_k))}{\partial x_i}+\frac{\partial s(\cosh(x_i),\cosh(x_l),\cosh(x_m))}{\partial x_i}.\]
Using the expression for the derivative of \(s\) (Lemma \ref{lem:hextrig}), we get
\begin{gather*}
\frac{\partial F_H}{\partial x_i}=\\
\sinh(x_i)\left(\cfrac{\cosh(x_i)-1-\cosh(x_j)-\cosh(x_k)}{(\cosh(x_i)-1)\sqrt{\cosh(x_i)^2+\cosh(x_j)^2+\cosh(x_k)^2+2\cosh(x_i)\cosh(x_j)\cosh(x_k)-1}}\right.\\+\left.\frac{\cosh(x_i)-1-\cosh(x_l)-\cosh(x_m)}{(\cosh(x_i)-1)\sqrt{\cosh(x_i)^2+\cosh(x_l)^2+\cosh(x_m)^2+2\cosh(x_i)\cosh(x_l)\cosh(x_m)-1}}\right)
\end{gather*}
So if \(x_i\leq x_j,x_k,x_l,x_m\) the derivative is negative. In particular, this holds if  \(\ell_{\alpha_i}(X)=\min_{j}\ell_X(\alpha_j)\), proving (1).

If \(\frac{\partial F_H}{\partial x_j}(\varphi_H(X))=0\), then \(\ell_{\alpha_i}(X)\) is a zero of the derivative of \(F_{H,X,i}\). Denote by \(c_j,c_k,c_l,c_m\) the hyperbolic cosines of \(\ell_{\alpha_j}(X),\ell_{\alpha_k}(X),\ell_{\alpha_l}(X),\ell_{\alpha_m}(X)\). A point \(y\in(0,\infty)\) is a zero of the derivative of \(F_{H,X,i}\) if and only if
\begin{gather*}
\frac{\cosh(y)-1-c_j-c_k}{(\cosh(y)-1)\sqrt{\cosh(y)^2+c_j^2+c_k^2+2\cosh(y)c_jc_k-1}}+\\\frac{\cosh(y)-1-c_l-c_m}{(\cosh(y)-1)\sqrt{\cosh(y)^2+c_l^2+c_m^2+2\cosh(y)c_lc_m-1}}=0
\end{gather*}
which implies
\[
\def\arraystretch{2.2}
\left.\begin{array}{l}\displaystyle{(\cosh(y)-1-c_j-c_k)\sqrt{\cosh(y)^2+c_l^2+c_m^2+2\cosh(y)c_lc_m-1}}=\\-(\cosh(y)-1-c_l-c_m)\sqrt{\cosh(y)^2+c_j^2+c_k^2+2\cosh(y)c_jc_k-1}\end{array}\right.
\]
which in turn implies
\begin{gather*}
(\cosh(y)-1-c_j-c_k)^2(\cosh(y)^2+c_l^2+c_m^2+2\cosh(y)c_lc_m-1)=\\(\cosh(y)-1-c_l-c_m)^2(\cosh(y)^2+c_j^2+c_k^2+2\cosh(y)c_jc_k-1).
\end{gather*}
This is a third degree equation in \(z=\cosh(y)\) (the coefficients of \(\cosh(y)^4\) cancel out), which has at most three real solutions, one of which is \(z=1\). Since \(z=1\) implies that \(y=0\), which is impossible, there are at most two zeroes of the derivative of \(F_{H,X,i}\) in the interval \((0,\infty)\). Since by Lemma \ref{lem:hextrig}
\[\lim_{y\to 0}F_{H,X,i}(y)=+\infty=\lim_{y\to +\infty}F_{H,X,i}(y),\]
there should be at least one zero of the derivative at the global minimum, and if there are two the other zero can only be an inflexion point.
\end{proof}

\begin{prop}\label{prop:localmax}
Let \(L>0\) and \(S\) be a surface of signature \((g,n)\). If \(X\in\teich(S;L)\) is a local maximum for the orthosystole, then \(\okiss(X)=6g-6+3n\) (ie \(\okiss(X)\) is maximal).
\end{prop}

\begin{proof}

Let \(X\in\teich(S;L)\) be a surface with \(\okiss(X)=k<6g-6+3n\). Let \(M=\{\alpha_1,\dots,\alpha_k\}\) be the collection of orthosystoles of \(X\). Complete \(M\) to a hexagon decomposition \(H=\{\alpha_1,\dots,\alpha_{6g-6+3n}\}\). Our goal is to show that \(X\) is not a local maximum.

The idea of the proof is to define surfaces arbitrarily close to \(X\) by increasing the lengths of all arcs in \(M\) and varying the length of \(\alpha_{6g-6+3n}\) in such a way that the boundary still has total length \(L\). We want to vary all these lengths by a sufficiently small amount so that, by Lemma \ref{lem:finitelymanylengths}, the orthosystole of the new surfaces is the minimum of the lengths of the curves in \(M\), which in turn implies that the surfaces we construct all have orthosystole strictly longer than \(\osys(X)\). For this, we need to know that if we keep the total boundary length fixed and we vary continuously all lengths of arcs in \(H\) but one, the length of the last arc varies continuously in terms of the other lengths. Since the total boundary length of a surface is defined by an explicit function of the lengths of the arcs in \(H\) (the function \(F_H\) from Section \ref{sec:hexdec}), we would want to use the implicit function theorem to show this continuity property. For this, we need the derivative \(\frac{\partial F_H}{\partial x_{6g-6+3n}}(\varphi_H(X))\) to be nonzero, which is unfortunately not always the case. So when this derivative is zero, we will instead increase the lengths of \(\alpha_2,\dots,\alpha_k\) and \(\alpha_{6g-6+3n}\) by small controlled amounts and prove that to keep the total boundary length fixed we need to increase the length of \(\alpha_1\) as well.

Let us formalize the argument. For ease of notation, in what follows we set \(N:=6g-6+3n\). Note first that by Lemma \ref{lem:finitelymanylengths} and Theorem \ref{thm:parametrization} there is \(\varepsilon_1>0\) such that if \(Y\in\teich(S)\) and \[\max_{i\in\{1,\dots,N\}}|\ell_{\alpha_i}(Y)-\ell_{\alpha_i}(X)|<\varepsilon_1,\] then
\[\osys(Y)=\min_{i\in\{1,\dots,k\}}\ell_{\alpha_i}(Y).\]
We assume that \(\varepsilon_1\) is small enough so that for every \(i\), \(\ell_{\alpha_i}(X)-\varepsilon_1>0\).

Let \(a=(a_1,\dots,a_N):=\varphi_H(X)\).

\noindent\ul{Case 1:} \(\frac{\partial F_H}{\partial x_N}(a)\neq 0\).

\noindent By the implicit function theorem there is \(\varepsilon_2>0\) and a continuous function
\[h:U=\left\{(y_1,\dots,y_{N-1})\in\R^{N-1}\;\middle|\; \max_{1\leq i\leq N-1}|y_i-a_i|<\varepsilon_2\right\}\to \R\] such that \(F_H(y_1,\dots, y_{N-1},h(y_1,\dots,h_{N-1}))=L\). We can assume that \(\varepsilon_2\) is small enough so that \(U\subset (\R_+)^{N-1}\) and \(h(U)\subset (a_N-\varepsilon_1,a_N+\varepsilon_1)\).

For every \(\varepsilon<\min\{\varepsilon_1,\varepsilon_2\}\) define 
\[X_\varepsilon=\varphi_H^{-1}(a_1+\varepsilon,\dots,a_k+\varepsilon,a_{k+1},\dots,a_{N-1}, h(a_1+\varepsilon,\dots,a_k+\varepsilon,a_{k+1},\dots,a_{N-1})).\]
In other words, \(X_\varepsilon\) is the surface where the lengths of the orthogeodesics in \(M\) are all increased by \(\varepsilon\) and the length of \(\alpha_N\) is adjusted so that the total length of the boundary of \(X_\varepsilon\) is \(L\). By the choice of \(\varepsilon_1,\varepsilon_2\) and \(\varepsilon\)
\[\osys(X_\varepsilon)=\min_{1\leq i\leq k}\ell_{\alpha_i}(Y)=\min_{1\leq i\leq k}a_i+\varepsilon=\osys(X)+\varepsilon.\]
Since \(X_\varepsilon\) converges to \(X\) when 
\(\varepsilon\) goes to zero, \(X\) is not a local maximum of \(\osys\).

\noindent\ul{Case 2:} \(\frac{\partial F_H}{\partial x_{N}}(a)=0\).

\noindent By Lemma \ref{lem:derivatives}, \(\frac{\partial F_H}{\partial x_{i}}(a)<0\) for every \(1\leq i\leq k\) and \(a_{N}\) is a global minimum or an inflexion point of the function
\[y\mapsto F_H(a_1,\dots,a_{N-1},y).\]
In case it's an inflexion point, let us assume that the function is locally monotone increasing around \(a_{N}\); the case in which the function is locally decreasing is analogous.

By the implicit function theorem, there is \(\varepsilon_2,\varepsilon_3>0\) and a continuous function
\[h:U=\left\{(y_2,\dots,y_{N})\in\R^{N-1}\;\middle|\; \max_{2\leq i\leq N}|y_i-a_i|<\varepsilon_2\right\}\to \R\] such that \(F_H(h(y_2,\dots, y_{N}),y_2,\dots, y_N)=L\). We can assume that \(\varepsilon_2,\varepsilon_3\) are small enough so that \(U\subset (\R_+)^{N-1}\), \(h(U)\subset (a_1-\varepsilon_3,a_1+\varepsilon_3)\), for some \(\varepsilon_3<\varepsilon_1\), and
\[\frac{\partial F_H}{\partial x_1}(z,y_2,\dots,y_{N})<0\]
for all \((y_2,\dots,y_N)\in U\) and \(z\in(a_1-\varepsilon_3,a_1+\varepsilon_3)\).

For every \(\varepsilon>0\), \(\varepsilon<\min\{\varepsilon_1,\varepsilon_2\}\) there is \(\delta>0\), \(\delta<\min\{\varepsilon_1,\varepsilon_2\}\), such that the total length of the boundary of the surface obtained by increasing\footnote{If instead the function \(y\mapsto F_H(a_1,\dots,a_{N-1},y)\) is locally monotone decreasing around \(a_N\), we decrease the length of \(\alpha_N\) by \(\varepsilon\).} the length of \(\alpha_N\) by \(\varepsilon\) and the lengths of \(\alpha_2,\dots,\alpha_k\) by \(\delta\) is more than \(L\). This holds since the boundary length is a continuous function of the lengths of arcs in \(H\) and by Lemma \ref{lem:derivatives}, if we increase the length of \(\alpha_N\), the boundary length increases, while if we increase the length of any arc in \(M\), the boundary length decreases. Now let \(X_{\varepsilon,\delta}\) be the surface such that:
\[\left\{\begin{array}{l}
\ell_{\alpha_i}(X_{\varepsilon,\delta})=a_i+\delta \hspace{3cm} \forall\, i\in\{2,\dots, k\}\\
\ell_{\alpha_i}(X_{\varepsilon,\delta})=a_i \hspace{3,6cm} \forall\,  i\in\{k+1,\dots, N-1\}\\
\ell_{\alpha_{N}}(X_{\varepsilon,\delta})=a_{N}+\varepsilon\\
\ell_{\alpha_1}(X_{\varepsilon,\delta})=h(a_2+\delta,\dots,a_k+\delta,a_{k+1},\dots,a_{N-1},a_{N}+\varepsilon)
\end{array}\right.\]
that is, 
\begin{gather*}
\varphi_H\left(X_{\varepsilon,\delta}\right)=\\
(h(a_2+\delta,\dots,a_k+\delta,a_{k+1},\dots,a_{N-1},a_N+\varepsilon),
a_2+\delta,\dots,a_k+\delta,a_{k+1},\dots,a_{N-1},a_{N}+\varepsilon).
\end{gather*}
Note that \(h(a_2+\delta,\dots,a_k+\delta,a_{k+1},\dots,a_{N-1},a_{N}+\varepsilon)>a_1\) because \[F_H(a_1,a_2+\delta,\dots,a_k+\delta,a_{k+1},\dots,a_{N-1},a_{N}+\varepsilon)>L\] and by the assumption on \(\varepsilon_2\) (and its consequence on the derivative with respect to \(x_1\)) we need to increase \(a_1\) to decrease \(L\). So
\[\osys(X_{\varepsilon,\delta})=\min_{1\leq i\leq k}\ell_{\alpha_i}(X_{\varepsilon,\delta})>\osys(X).\]
As \(\varepsilon\) tends to zero, so does \(\delta\) and thus \(X_{\varepsilon,\delta}\) converges to \(X\), showing again that \(X\) is not a local maximum.
\end{proof}

We can now prove the equivalence of the statements in Theorem \ref{thm:n=1}:

\begin{proof}[Proof of Theorem \ref{thm:n=1}]
Recall that the four conditions to be proven equivalent are:
\begin{enumerate}
\item \(X\) is a global maximum for the orthosystole function on \(\ms(S;L)\);
\item \(X\) is a local maximum for the orthosystole function on \(\ms(S;L)\);
\item \(\okiss(X)=6g-6+3n\);
\item \(\displaystyle\osys(X)=2\sinh^{-1}\left(\frac{1}{2\sinh\left(\frac{L}{24g-24+12n}\right)}\right)\).
\end{enumerate}
It is clear that (1) implies (2). Proposition \ref{prop:localmax} shows that (2) implies (3) and Lemma \ref{lem:hexdecsamelength} shows that (3) implies (4).

To prove that (4) implies (1), we note that all local maxima have the same orthosystole length (since (2) implies (4)) and thus, as we know that \(\osys\) admits a global maximum (Proposition \ref{prop:localmax}), all local maxima are also global maxima. Moreover, this also shows that
\[\sup_{X\in\ms(S;L)}\osys(X)=2\sinh^{-1}\left(\frac{1}{2\sinh\left(\frac{L}{24g-24+12n}\right)}\right)\]
hence (4) implies (1).
\end{proof}

Note that the equivalence of (1), (3) and (4) follow from Bavard's work, though we've given here alternative proofs. The following Corollary will finish the proof of Theorem \ref{thm:n=1}.

\begin{cor}\label{cor:numberofmax}
For every surface \(S\) of signature \((g,1)\), the function \(\osys:\ms(S; \ell)\to\R\) has exactly
\[\frac{2 (6g-5)!}{12^g g!(3g-3)!}\]
local maxima (which are all global maxima).
\end{cor}

\begin{proof}
This follows from Theorem \ref{thm:n=1}, Remark \ref{rmk:hexdec}, the fact that homeomorphism classes of hexagon decomposition are in correspondence with homeomorphism classes of ideal triangulations of a surface with one puncture, and the computation of the number of homeomorphism classes of such triangulations from \cite{bv_counting}.
\end{proof}

\section{Surfaces with given boundary lengths}\label{sec:multiplecomponents}

Understanding local and global maxima of the orthosystole function in the case of multiple boundary components, each of whose length is fixed, seems significantly harder. The main reason is that we cannot in general guarantee that there is a surface with a hexagon decomposition of orthosystoles. Indeed, if a surface \(X\) of signature \((g,n)\) has a hexagon decomposition of orthosystoles, then each boundary component is a union of arcs of length \(\alpha\), where \(\alpha\) is the length of any other side of a right-angled hexagon of alternating sides of length \(\osys(X)\). In particular, if the boundary lengths are not rational multiples of each other there won't be any hexagon decomposition of orthosystoles. But even when the boundary lengths are rational multiples of each other it is not clear if it is always possible to have a hexagon decomposition of orthosystoles.

As mentioned before, in \cite{bavard_anneaux} it is the sum of the boundary lengths to be fixed, instead of the individual lengths. Note that for every signature and every sum of boundary lengths we can find a hyperbolic surface with a hexagon decomposition of orthosystoles, which means that it attains the bound in Theorem \ref{thm:Bavard}.

One case in which we can show the existence of such a hexagon decomposition is the case of surfaces with two boundary components of the same length. We will then use this construction to construct surfaces with large orthosystole and prove Theorem \ref{thm:natleast2}.

\begin{lemma}\label{lem:n=2samelength}
Let \(S_{g,2}\) be a surface of signature \((g,2)\). For any \(\ell>0\) there is \(X\in\ms(S_{g,2};\ell,\ell)\) with a hexagon decomposition of orthosystoles and 
\[\osys(X)=\cosh^{-1}\left(\frac{\cosh\left(\frac{\ell}{6g}\right)}{\cosh\left(\frac{\ell}{6g}\right)-1}\right)=\max_{Y\in\ms(S_{g,2};\ell,\ell)}\osys(Y).\]
\end{lemma}

\begin{proof}
Color one boundary component of \(S_{g,2}\) blue and the other red. We will first show, by induction on the genus, that \(S_{g,2}\) has a hexagon decomposition such that half of the hexagons have two red sides and one blue side and half of the hexagons have two blue sides and one red side.

In the base case, \(g=1\), an example of such a hexagon decomposition is given in Figure \ref{fig:2red1bluegenus1}.

\begin{figure}[h]
\begin{center}
\includegraphics{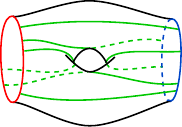}
\caption{A hexagon decomposition on \(S_{1,2}\)}\label{fig:2red1bluegenus1}
\end{center}
\end{figure}

For the induction step, assume \(H\) is a hexagon decomposition with the required properties on \(S_{g,2}\). Choose a hexagon \(H_1\) with two red sides and a hexagon \(H_2\) with two blue sides. Remove a disk from the interior of \(H_1\) and one from the interior of \(H_2\) and glue the two new boundary components together. We now have a surface of signature \((g+1,2)\) with a collection of disjoint arcs. We add arcs as in Figure \ref{fig:2red1blueg+1} to get a hexagon decomposition with the required properties.

\begin{figure}[h]
\begin{center}
\includegraphics{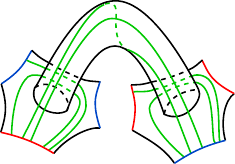}
\caption{From genus \(g\) to genus \(g+1\)}\label{fig:2red1blueg+1}
\end{center}
\end{figure}

So for every \(g\geq 1\) we can choose a hexagon decomposition with the properties above; if we choose all arcs in this hexagon decomposition to have the same length \(a\) we get a hyperbolic surface \(X\) with a hexagon decomposition of orthosystoles (Lemma \ref{lem:decosys}), whose boundary components have the same length, equal to \(6g\alpha\), where \(\alpha\) is the length of the other sides of a hexagon with alternating sides of length \(a,a,a\). So the boundary lengths are \(\ell\) if and only if
\[\ell=6g\cosh^{-1}\left(\frac{\cosh(a)}{\cosh(a)-1}\right)\]
and hence it suffices to choose \(a\) to be the unique solution of:
\[\cosh(a)=\frac{\cosh\left(\frac{\ell}{6g}\right)}{\cosh\left(\frac{\ell}{6g}\right)-1}.\]
By Lemma \ref{lem:decosys}
\[\osys(X)=\cosh^{-1}\left(\frac{\cosh\left(\frac{\ell}{6g}\right)}{\cosh\left(\frac{\ell}{6g}\right)-1}\right).\]
By Theorem \ref{thm:Bavard}, \(X\) is a global maximum for \(\osys\) on \(\ms(S_{g,2};2\ell)\), and thus on \(\ms(S_{g,2};\ell,\ell)\).
\end{proof}

Using Lemma \ref{lem:n=2samelength}, we can prove part (2) of Theorem \ref{thm:natleast2}:

\begin{proof}[Proof of Theorem \ref{thm:natleast2}, part (2)]
Suppose first \(n=2\) and the genus is odd. Let \(X_1\) and \(X_2\) be hyperbolic surfaces of signature \(\left(\frac{g-1}{2},2\right)\) and boundary components \(\gamma_1,\delta_1\) and \(\gamma_2,\delta_2\) respectively, such that:
\begin{itemize}
\item \(\ell_{\gamma_i}(X_i)=\ell_{\delta_i}(X_i)=\ell_i\), and
\item \(X_i\) has a hexagon decomposition of orthogeodesics of minimal length as in Lemma \ref{lem:n=2samelength}.
\end{itemize} Let \(Y\) be a two-holed torus, with boundary components \(\eta_1\) and \(\eta_2\) of lengths \(\ell_1\) and \(\ell_2\) respectively. Let \(X\) be a surface obtained from the \(X_i\)'s and \(Y\) by gluing, for every \(i\), \(\delta_i\) to \(\eta_i\).
 Note that any orthogeodesic \(\alpha\) in \(X\) is either contained in a \(X_i\) or it contains a subarc joining two boundary components of a \(X_i\). By construction each \(X_i\) contains an orthogeodesic from \(\gamma_i\) to itself of length \(\osys(X_i)\), which gives an orthogeodesic on \(X\) of the same length. Thus
\begin{align*}
\osys(X)&=\min\{\osys(X_1),\osys(X_2)\}=\osys(X_2)=\\
&=\cosh^{-1}\left(\frac{\cosh\left(\frac{\ell_2}{3g-3}\right)}{\cosh\left(\frac{\ell_2}{3g-3}\right)-1}\right)=\cosh^{-1}\left(\cfrac{\cosh\left(\cfrac{\ell_2}{6\left\lfloor \frac{g}{2}\right\rfloor}\right)}{\cosh\left(\cfrac{\ell_2}{6\left\lfloor \frac{g}{2}\right\rfloor}\right)-1}\right),
\end{align*}
where the second equality follows from the monotonicity of
\[x\mapsto \cosh^{-1}\left(\frac{\cosh(x)}{\cosh(x)-1}\right)\]
for \(x>0\).

The case \(n=2\) and even genus is similar; the only difference is that we choose \(S_1\) to have signature \(\left(\frac{g}{2},2\right)\) and \(S_2\) to have signature \(\left(\frac{g}{2}-1,2\right)\). We then get a surface with orthosystole
\[\cosh^{-1}\left(\frac{\cosh\left(\frac{\ell_2}{3g-6}\right)}{\cosh\left(\frac{\ell_2}{3g-6}\right)-1}\right)\geq \cosh^{-1}\left(\cfrac{\cosh\left(\cfrac{\ell_2}{6\left\lfloor \frac{g}{2}\right\rfloor}\right)}{\cosh\left(\cfrac{\ell_2}{6\left\lfloor \frac{g}{2}\right\rfloor}\right)-1}\right),\]
where again we are using the motononicity of
\[x\mapsto \cosh^{-1}\left(\frac{\cosh(x)}{\cosh(x)-1}\right)\]
for \(x>0\).

Suppose now \(n\geq 3\). Let \(q=\lfloor \frac{g}{n}\rfloor\) and let \(r\) be such that \(g=qn+r\). Consider:
\begin{itemize}
\item \(X_i\) a hyperbolic surface of signature \(q,2\), with boundary components \(\gamma_i,\delta_i\) of length \(\ell_i\) and with a hexagon decomposition as in Lemma \ref{lem:n=2samelength};
\item \(Y\) a hyperbolic surface of signature \(q,n\) and boundary components \(\eta_1,\dots\eta_n\) of lengths \(\ell_1,\dots,\ell_n\) respectively.
\end{itemize}
Let \(X\) be a surface obtained from the \(X_i 's\) to \(Y\) by gluing, for every \(i\), \(\delta_i\) to \(\eta_i\) with any choice of twist. Note that \(X\in\teich(S;\ell_1,\dots,\ell_n)\); furthermore, the same argument as above shows that 
\[\osys(X)=\osys(Y_n)=\cosh^{-1}\left(\cfrac{\cosh\left(\cfrac{\ell_n}{6\left\lfloor \frac{g}{n}\right\rfloor}\right)}{\cosh\left(\cfrac{\ell_n}{6\left\lfloor \frac{g}{n}\right\rfloor}\right)-1}\right).\]
\end{proof}

With similar techniques to those used in the proof of Proposition \ref{prop:localmax}, we can show that the surfaces constructed in the previous proof aren't even local maxima for the orthosystole function. 

On the other hand, there are choices of signatures and boundary lengths so that, in the corresponding moduli space, we can construct surfaces with a hexagon decomposition of orthosystoles, and thus maximizing the orthosystole function. The following lemma gives an explicit example of such a construction.

\begin{lemma}\label{lem:infinitefamily}
Fix \(n\geq 1\); for any \(m\geq 0\) let \(S_m\) be a surface of signature \((nm+1,n)\). Then for every \(\ell>0\) there is \(X_m\in\ms(S_m;\ell,\dots,\ell)\) with a hexagon decomposition of orthogeodesics of minimal length. In particular,
\(X_m\) maximizes \(\osys\) in \(\ms(S_m;\ell,\dots,\ell)\)
\[\osys(X_m)=\cosh^{-1}\left(\frac{\cosh\left(\frac{\ell}{12m+6}\right)}{\cosh\left(\frac{\ell}{12m+6}\right)-1}\right).\]
\end{lemma}

\begin{proof}
Fix \(m\geq 0\); the choice of genus allows us to find an order \(n\) symmetry \(\varphi\) of \(S_m\), such that a fundamental domain for the action is a subsurface \(F\) of signature \((m,3)\), as in Figure \ref{fig:symmetry}.

\begin{figure}[h]
\begin{center}
\includegraphics{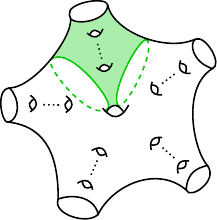}
\caption{The case \(n=5\); the fundamental domain is shaded.}\label{fig:symmetry}
\end{center}
\end{figure}
We then construct a hexagon decomposition of \(S_m\) as follows.

Look first at \(F\) and denote by \(\gamma\) its boundary component which is also a boundary component of \(S_m\). Let \(\hat{F}\) be the surface obtained by gluing two disks to the two other boundary components of \(F\). Fix a hexagon decomposition of \(\hat{F}\) and homotope the arcs so that they are contained in \(F\) and that the two disks in \(\hat{F}\ssm F\) are contained in different hexagons. Look now at \(F\); for each hexagon containing a boundary component add two new essential arcs from \(\gamma\) to itself, as in Figure \ref{fig:adding-arcs}.

\begin{figure}[]
\begin{center}
\includegraphics{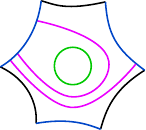}
\caption{The sides in blue are subarcs of \(\gamma\); the arcs in purple are to be added.}\label{fig:adding-arcs}
\end{center}
\end{figure}

We now have a collection of \(6m-3+4=6m+1\) arcs on \(F\). Take all the arcs in their \(\varphi\)-orbit, so that we have \(6mn+n\) arcs on \(S_m\). It is not hard to see that there are \(n\) complementary components \(C_i\), for \(i=0,\dots, n-1\), which are not hexagons, where \(C_i\) contains \(\varphi^{i}(\gamma)\) and \(\varphi^{i+1}(\gamma)\). For every \(i\), add two arcs in \(C_i\), each joining \(\varphi^{i}(\gamma)\) to \(\varphi^{i+1}(\gamma)\). The union of all these arcs is a hexagon decomposition \(H_m\) of \(S_m\) with the property that each boundary component of \(S_m\) is split into the same number of subarcs. This implies that if all arcs in \(H_m\) are given the same length, all boundary components of \(S_m\) will have the same length. To conclude the proof we just need to show that for every \(\ell>0\) we can choose a length \(a=a(\ell,n,m)>0\) so that if all arcs in \(H_m\) have length \(a\), each boundary component of \(S_m\) has length \(\ell\). Since \(H_m\) contains \(4nm+2n\) hexagons, by Lemma \ref{lem:hextrig} we know that \(a\) needs to satisfy
\[n\ell=(4nm+2n)s(\cosh(a),\cosh(a),\cosh(a))\]
ie
\[n\ell=(12nm+6n)\cosh^{-1}\left(\frac{\cosh(a)}{\cosh(a)-1}\right)\]
which implies
\[\cosh(a)=\frac{\cosh\left(\frac{\ell}{12m+6}\right)}{\cosh\left(\frac{\ell}{12m+6}\right)-1}.\]
As the right-hand side is bigger than one, there is a unique solution. Then \(X_m\) can be chosen to be the surface where the arcs in \(H_m\) have length \(a\) and by Lemma \ref{lem:decosys} they are the orthosystoles.

By Theorem \ref{thm:Bavard}, these surfaces are global maximizers for the orthosystole in \(\ms(S_m;n\ell)\), and thus in \(\ms(S_m;\ell,\dots,\ell)\).
\end{proof}

We end this section by using Lemma \ref{lem:n=2samelength} to show that there is no version of Mumford's compactness theorem that can be stated using orthosystoles instead of systoles:

\begin{lemma}\label{lem:noMumford}
Let \(\ell>0\). For every \(R>0\) there is \(g_0\geq 1\) so that for every \(g\geq g_0\), if \(S\) is a surface of signature \((g,1)\), the set
\[\{X\in\ms(S;\ell)\st \osys(X)\geq R\}\]
is not compact.
\end{lemma}

\begin{proof}
By Lemma \ref{lem:n=2samelength}, there is a surface \(X_g\) of signature \((g-1,2)\) and boundary lengths \(\ell,\ell\) with a hexagon decomposition of orthosytoles and 
\[\osys(X_g)=\cosh^{-1}\left(\frac{\cosh\left(\frac{\ell}{6(g-1)}\right)}{\cosh\left(\frac{\ell}{6(g-1)}\right)-1}\right).\]
For any one-holed torus \(T\) of boundary length \(\ell\), we can construct a surface \(Y_{g,T}\) by gluing \(T\) to a boundary component of \(X_g\) with any choice of twist parameter. Note that
\[\osys(Y_{g,T})\geq \osys(X_g),\]
because any orthogeodesic in \(Y_{g,T}\) is either an orthogeodesic of \(X_g\) or it contains an arc joining two boundary components of \(X_g\).

Since \(\osys(Y_{g,T})\to\infty\) as \(g\to\infty\), for every \(R>0\) there is \(g_0\) so that \(\osys(Y_{g,T})>R\) for every \(g\geq g_0\). It is sufficient to choose a sequence \(T_n\) of one-holed tori with a curve shrinking to zero to get a sequence \(Y_{g,T}\) leaving every compact.
\end{proof}

\section{Orthosystole maximizers and orthokissing number}\label{sec:osys&okiss}

In this section we prove part (3) of Theorem 
\ref{thm:natleast2}; that is, 
a bound on the orthokissing number of a surface which is a local maximizer for the orthosystole. The main result we need is the following:

\begin{prop}\label{prop:filling}
Let \(S\) be a surface of signature \((g,n)\) and let \(0<\ell_1,\dots,\ell_n\). If \(X\in\teich(S;\ell_1,\dots, \ell_n)\) is a local maximum for the orthosystole function, the collection of orthosystoles of \(X\) fill the surface, ie every simple closed geodesic on \(X\) intersects at least one orthosystole.
\end{prop}
\begin{proof}
By contradiction, suppose not. Denote by \(\gamma_1,\dots,\gamma_n\) the boundary components of \(S\). Let \(\alpha\) be a simple closed geodesic disjoint from all orthosystoles of \(X\). Cut \(X\) along \(\alpha\) to get a (possibly disconnected) surface \(Y\) with \(n+2\) boundary components, \(n\) corresponding to \(\gamma_1,\dots,\gamma_n\) and two new ones, which we denote by \(\alpha_1\) and \(\alpha_2\). By Theorem \ref{thm:lengthexpansion}, for every \(\varepsilon>0\) there is a surface \(Y_\varepsilon\) satisfying
\[\ell_{\gamma_i}(Y_\varepsilon)=\ell_{\gamma_i}(Y)\;\;\forall i=1,\dots, n\]
and
\[\ell_{\alpha_i}(Y_\varepsilon)=\ell_{\alpha_i}(Y)+\varepsilon\;\;\forall i=1,2\]
and such that, by Remark \ref{rmk:lengthexpansion} and since orthosystoles of \(X\) are disjoint from \(\alpha\), for any orthosystole \(\beta\) of \(X\)
\[\ell_\beta(Y_{\varepsilon})>\ell_\beta(Y).\]
Glue back \(Y_\varepsilon\), with the same twist parameter about \(\alpha\) as \(X\), to get a surface \(X_\varepsilon\) such that for any orthosystole \(\beta\) of \(X\)
\[\ell_\beta(X_\varepsilon)>\ell_\beta(X).\]
For \(\varepsilon\) small enough, the set of orthosystoles of \(X_\varepsilon\) is a subset of the set of orthosystoles of \(X\) (by discreteness of the orthogeodesic spectrum), for every \(\varepsilon\) small enough we get
\[\osys(X_\varepsilon)>\osys(X),\]
contradicting local maximality of \(X\).
\end{proof}

A consequence of the proposition is that local maximizers of the orthosystole functions lie in a compact subset of moduli space.

\begin{cor}
Let \(S\) be a surface of signature \((g,n)\) and let \(0\leq \ell_1\leq \dots\leq \ell_n\). Then there is \(\varepsilon>0\) such that any local maximum for the orthosystole function on \(\ms(S;\ell_1,\dots,\ell_n)\) lies in \(\ms_\varepsilon(S;\ell_1,\dots,\ell_n)\). We can choose
\[\varepsilon=2\sinh^{-1}\left(2\sinh\left(\frac{\ell(\partial X)}{24g-24+12n}\right)\right).\]
\end{cor}

\begin{proof}
If \(X\) is a local maximum, its orthosystoles fill. So the systole of \(X\) intersects an orthosystole of \(X\) and thus
\[\osys(X)\geq 2w(\sys(X))=2\sinh^{-1}\left(\frac{1}{\sinh(\sys(X)/2)}\right).\]
By Theorem \ref{thm:Bavard}
\[\osys(X)\leq 2\sinh^{-1}\left(\frac{1}{2\sinh\left(\frac{\ell(\partial X)}{24g-24+12n}\right)}\right).\]
Combining the two inequalities yields the result.
\end{proof}

To deduce the lower bound on the orthokissing number of a local maximizer for the orthosystole function from Proposition \ref{prop:filling} we just need to compute how many disjoint arcs are needed to fill:

\begin{lemma}\label{lem:filling}
Let \(S\) be a surface of signature \((g,n)\). Let \(\A\) be a filling collection of disjoint arcs of \(S\). Then
\[|\A|\geq \left\{
\begin{array}{ll}
2g-2+n & \text{if } n\geq 2\\
2g & \text{if } n=1
\end{array}
\right.
\]
and the bound is sharp.
\end{lemma}

\begin{proof}
Suppose \(|\A|=k\). As \(\A\) is filling, it cuts \(S\) into \(N_1\) topological disks and \(N_2\) peripheral annuli. Denote the disks by \(D_1,\dots, D_{N_1}\) and the annuli by \(A_1,\dots, A_{N_2}\).

Each disk \(D_i\) is a polygon with sides alternating between arcs of \(\A\) and arcs of the boundary. Let \(d_i\) be the number of arcs of \(\A\), so that \(D_i\) is a \(2d_i\)-gon.

One boundary component of an annulus \(A_i\) is a boundary component of \(S\) and the other is a polygon with sides alternating between arcs of \(\A\) and arcs of the boundary. Let \(a_i\) be the number of arcs on \(\A\) in the boundary of \(A_i\). Pick an arc \(\beta_i\) from one boundary component of \(A_i\) to the other, so that it intersects the polygonal boundary component of \(A_i\) in the middle of an arc of the boundary of \(S\) (see Figure \ref{fig:disks_annuli}). Cutting \(A_i\) along \(\beta_i\) yields a polygon \(B_i\) with \(2(a_i+2)\) sides, half of which are arcs of the boundary of \(S\).

\begin{figure}[h]
\begin{center}
\begin{overpic}{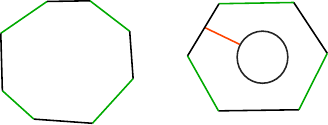}
\end{overpic}
\caption{A disk and an annulus in \(S\ssm\A\). The green sides are arcs of \(A\), the black ones are arcs of the boundary of \(S\) and the orange one is a \(\beta_i\).}\label{fig:disks_annuli}
\end{center}
\end{figure}

We now compute the Euler characteristic using the CW-complex decomposition of \(S\) given by the \(D_i\) and the \(B_i\). Note that each arc of \(A\) appears in exactly two polygons, so
\[\sum_{i=1}^{N_1}d_i+\sum_{i=1}^{N_2}a_i=2k.\]

Each vertex of a \(D_i\) or \(B_i\) is identified with another vertex, so the number of \(0\)-cells is
\[\frac{1}{2}\left(\sum_{i=1}^{N_1}(2d_i)+\sum_{i=1}^{N_2}(2(a_i+2)\right)=\sum_{i=1}^{N_1}d_i+\sum_{i=1}^{N_2}a_i+\sum_{i=1}^{N_2}2=2k+2N_2.\]

The edges of \(D_i\) or \(B_i\) not coming from arcs of the boundary of \(S\) are identified in pairs, while the other edges are not identified with any other edge. Thus the number of \(1\)-cells is
\[\frac{1}{2}\left(\sum_{i=1}^{N_1}d_i+\sum_{i=1}^{N_2}(a_i+2)\right)+\sum_{i=1}^{N_1}d_i+\sum_{i=1}^{N_2}(a_i+2)=k+N_2+2k+2N_2=3k+3N_2.\]

Finally, the number of \(2\)-cells is \(N_1+N_2\). So
\[2-2g-n=2k+2N_2-(3k+3N_2)+N_1+N_2=-k+N_1,\]
which implies that \(k\geq 2g-2+n+N_1\).

If \(n=1\), all arcs start and end at the unique boundary component, so there cannot be annuli in the complement of \(\A\), ie \(N_2=0\) and \(N_1\geq 1\). So \(k\geq 2g\), with equality if and only if \(N_1=1\).

If \(n\geq 2\), \(k\geq 2g-2+n\), with equality if and only if \(N_1=0\).

To prove that the two bounds are sharp, we just need to exhibit:
\begin{itemize}
\item if \(n=1\), a collection of arcs whose complement is a single polygon;
\item if \(n\geq 2\), a collection of arcs whose complement is a union of annuli.
\end{itemize}
These are given in Figure \ref{fig:minimalfilling}.

\begin{figure}[h]
\begin{center}
\begin{overpic}{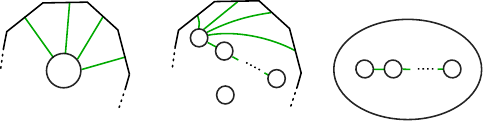}
\put(10,-5){\(n=1\)}
\put(-1,22){\(a_1\)}
\put(12,26){\(b_1\)}
\put(22,22){\(a_1^{-1}\)}
\put(26,15){\(b_1^{-1}\)}
\put(60,-5){\(n\geq 2\)}
\put(35,22){\(a_1\)}
\put(48,26){\(b_1\)}
\put(58,22){\(a_1^{-1}\)}
\put(61,15){\(b_1^{-1}\)}
\end{overpic}
\vspace{.2cm}
\caption{Filling collections of arcs of minimal size}\label{fig:minimalfilling}
\end{center}
\end{figure}
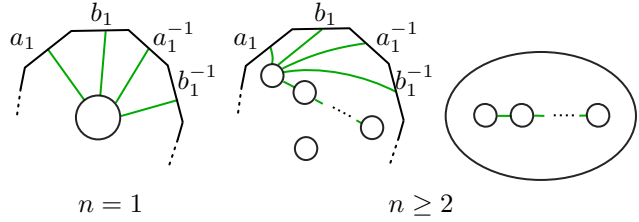
\end{proof}

Part (3) of Theorem \ref{thm:natleast2} is now an easy corollary of the results of this section:

\begin{cor}
Let \(S\) be a surface of signature \((g,n)\), for \(n\geq 2\), and let \(0<\ell_1\leq \dots\leq \ell_n\). Then if \(X\in\teich(S;\ell_1,\dots,\ell_n)\) is a local maximum for the orthosystole function, \(\okiss(S)\geq 2g-2+n\).
\end{cor}

\begin{proof}
This follows from Proposition \ref{prop:filling}, Lemma \ref{lem:simple&disjoint} and Lemma \ref{lem:filling}.
\end{proof}

\appendix \section{Orthosystole and injectivity radius of the boundary}
The goal of this appendix is to explain why Bavard's result in \cite{bavard_anneaux} can be stated in terms of orthosystoles and orthokissing numbers.

Let us start with the setup by Bavard. He considers a compact hyperbolic surface \(X\) and a collection \(c\) of pairwise disjoint simple closed geodesics, which is assumed to contain the boundary components of \(X\), if there are any. For such a collection, he defines \(\Sigma=\Sigma(c)\subset X\) to be the union of all points \(x\in X\) such that there are at least two (geodesic) paths from \(x\) to \(c\) realizing the distance of \(x\) from \(c\). He proves (\cite[Lemme 1]{bavard_anneaux}) that \(\Sigma\) is a graph whose edges are geodesic segments. The \emph{injectivity radius} of \(c\) is
\[\inj(c):=\min_{x\in\Sigma}d(x,c).\]

We recall here the statements from (1), (3) and (4) of \cite[Th\'eor\`eme 1]{bavard_anneaux} for the injectivity radius of \(c=\partial X\). For ease of notation, we simply write \(\Sigma\) instead of \(\Sigma(\partial X)\).

\begin{thm}
Let \(X\) be a compact hyperbolic surface with non-empty boundary of signature \(g,n\) and let $L$  be the sum of its boundary lengths. If \(r\) is the injectivity radius of \(\partial X\), then
\[\sinh(r)\sinh\left(\frac{L}{12(6g-6+3n)}\right)\leq\frac{1}{2},\]
with equality if and only if all edges of \(\Sigma\) have the same length and meet at \(\frac{2\pi}{3}\) angles. Moreover, for every $L$, there is a hyperbolic surface attaining the bound.
\end{thm}

The equivalence of this theorem and Theorem \ref{thm:Bavard} follows from the next two lemmas:

\begin{lemma}\label{lem:inj&osys}
If \(c=\partial X\), \(\inj(c)=\frac{\osys(X)}{2}\).
\end{lemma}
\begin{proof}
If \(\alpha\) is an orthosystole, the midpoint \(x\) of \(\alpha\) belongs to \(\Sigma\) and thus \(\frac{\osys(X)}{2}\geq \inj(c)\). Conversely, if \(x\in \Sigma\) realizes the injectivity radius of \(c\), look at two distance-realizing paths \(a_1\) and \(a_2\) from \(x\) to \(c\). The concatenation \(a_2\ast a_1^{-1}\) is an essential arc from the boundary to the boundary, so if \(\alpha\) is the orthogeodesic in the homotopy class, we have
\[\osys(X)\leq \ell_\alpha(X)\leq \ell_{a_1}(X)+\ell_{a_2}(X)=2 d(x,c)=2\inj(c).\]
\end{proof}

\begin{lemma}\label{lem:Sigma&okiss}
Suppose \(c=\partial X\); then \(\Sigma\) has edges of the same length which meet at \(\frac{2\pi}{3}\) angles if and only if \(\okiss(X)=6g-6+3n\).
\end{lemma}
\begin{proof}
Suppose first that \(\Sigma\) has edges of the same length which meet at \(\frac{2\pi}{3}\) angles. Drop the perpendiculars from the vertices of \(\Sigma\) to \(c\); by \cite[Lemme 1]{bavard_anneaux}, we get a decomposition of \(X\) into isometric quadrilateral with angles \(\frac{\pi}{2},\frac{\pi}{2},\frac{\pi}{3},\frac{\pi}{3}\), where the side between the right angles is a segment of a boundary component and the opposite one is an edge of \(\Sigma\). In particular, for every edge of \(\Sigma\), there is an orthogeodesic which intersects it once orthogonally at its midpoint. By construction of \(\Sigma\), this orthogeodesic has length \(2\inj(c)\), so by Lemma \ref{lem:inj&osys} it is an orthosystole. As observed after the statement of \cite[Th\'eor\`eme 1]{bavard_anneaux}, since \(X\) deformation retracts onto \(\Sigma\), \(2-2g-n=v-e\), where \(v\) is the number of vertices of \(\Sigma\) and \(e\) the number of edges. As \(\Sigma\) is trivalent, \(v=\frac{2}{3}e\). So there are at least
\[e=6g-6+3n\]
orthosystoles, and by Lemma \ref{lem:simple&disjoint} this means that \(\okiss(X)=6g-6+3n\).

Conversely, if \(\okiss(X)=6g-6+3n\), then \(X\) admits a hexagon decomposition of orthosystoles, and it is not hard to show that \(\Sigma\) is the graph obtained by gluing one tripod per hexagon, and the tripod is given by the barycenter and the three perpendiculars from the barycenter to the three orthosystoles. As all hexagons are isometric, \(\Sigma\) satisfies the conditions in the statement.
\end{proof}

\bibliographystyle{plain}
\bibliography{references}

\end{document}